\documentclass[a4paper,11pt]{article}
\usepackage[utf8]{inputenc}

\usepackage{amsmath}
\usepackage{amsfonts}
\usepackage{amssymb}

\usepackage[T1]{fontenc}
\usepackage{gensymb}
\usepackage{lmodern}

\usepackage[left=2.5cm,right=2.5cm,top=2.5cm,bottom=3.2cm]{geometry}
\usepackage{apacite} 
\usepackage{natbib}

\usepackage[hyphens]{url}
\usepackage{hyperref}
\usepackage[hyphenbreaks]{breakurl}

\author{
  Dennis M\"uller\footnote{RWTH Aachen University. \texttt{dennis.mueller3@rwth-aachen.de} .} \ \  Maurice Chiodo\footnote{King's College, University of Cambridge. \texttt{mcc56@cam.ac.uk} .} \ \ James Franklin\footnote{University of New South Wales, Sydney. \texttt{j.franklin@unsw.edu.au} .}
}

\title{A Hippocratic Oath for mathematicians? \\Mapping the landscape of ethics in mathematics }
\begin{document}

\maketitle
\abstract{   
While the consequences of mathematically-based software, algorithms and strategies have become ever wider and better appreciated, ethical reflection on mathematics has remained primitive. We review the somewhat disconnected suggestions of commentators in recent decades with a view to piecing together a coherent approach to ethics in mathematics. Calls for a Hippocratic Oath for mathematicians are examined and it is concluded that while lessons can be learned from the medical profession, the relation of mathematicians to those affected by their work is significantly different. There is something to be learned also from the codes of conduct of cognate but professionalised quantitative disciplines such as engineering and accountancy, as well as from legal principles bearing on professional work. We conclude with recommendations that professional societies in mathematics should sponsor an (international) code of ethics, institutional mission statements for mathematicians and syllabuses for incorporation into mathematics degrees.  }

\section*{Introduction}

\let\thefootnote\relax\footnotetext{2020 \textit{AMS Classification:} 01A80. 00A30.}
\let\thefootnote\relax\footnotetext{\textit{Keywords:} Ethics in Mathematics, Hippocratic Oaths, Code of Ethics, Mathematics and Society.}
\let\footnote\relax\footnotetext{Submitted to \textit{Science and Engineering Ethics}.}

The contrast is stark between the massively ramifying consequences of mathematical work and the paucity of ethical reflection on it. It has become almost a platitude that “algorithms are ruling our lives” in areas ranging from social media feeds to healthcare decisions to predicting parolee recidivism. Algorithms honed by mathematicians will drive cars and fight wars. The “weapons of math destruction” of O’Neil’s popular book (2016) are recognised by those outside mathematics to be a major ethical concern. Ethical evaluation of them needs to be informed by mathematical expertise, since uninformed criticism will miss the mark. The more ethically positive uses of mathematics, such as CT scanning, need also to be appreciated, and mathematical technologies with potentially both good and bad uses such as cryptography need to be evaluated from a perspective of sound technical knowledge combined with ethical insight.

Yet the mathematical community does not have a solid tradition of ethical reflection on itself and the effects of its work. Cognate quantitative disciplines like engineering and accountancy (and even a subfield of mathematics itself, the actuarial profession) have well-developed codes of ethics and enforcement by professional bodies which include ethics in rigorous processes of professional accreditation. The Hippocratic Oath in medicine is the best-known of professional codes, a matter of pride among doctors since ancient times as reflecting the high impact of their work on patients. 

Mathematicians have been, as we will describe below, very slow to follow the lead of those in other professions of high societal impact. Although professional societies have developed some codes of conduct for practices internal to the academic profession (such as warning against plagiarism) little has been done about the wider issues of mathematical products that have major effects on society.  Calls for wider ethical codes for mathematicians, especially applied mathematicians, and calls for a Hippocratic Oath for mathematicians, have been muted, rare and disconnected.

In this article, we survey in Part I the development of what ethical reflection there has been on mathematical applications. In Part II we consider whether a Hippocratic Oath is an appropriate model for mathematicians. In Part III we explain what can be learned from codes of ethics in cognate quantitative disciplines and from certain basic ethical principles that have been incorporated into the law governing products that put users at risk. In our conclusion, we outline specific steps for the mathematical community to take in order to go beyond codes of conduct and to establish an environment actively conducive for ethical mathematical practice. These include institutional mission statements, integrated ethics courses in the mathematical curriculum and an international code of ethics.

For the purposes of this paper, we will use the working definition that “(practicing) mathematicians” are individuals who have a degree or training related to mathematics, and who currently apply mathematical skills and training in their everyday work.

We begin however with a brief consideration of the pure/applied distinction in mathematics, which is central to the discipline and has sometimes impeded ethical reflection on mathematics.

\subsection*{Pure and applied mathematics: ethical differences}

The distinction between pure and applied mathematics is conceptually and institutionally significant. Plainly, it impacts on ethical issues, since mathematical work applied to extra-mathematical problems has a far greater impact on society than research on purely mathematical questions.

To complicate matters, the meaning of “pure” and “applied” to mathematicians themselves is not the same as the naïve meanings of those words. Mathematics is institutionally divided into pure and applied primarily by content, rather than by what it is used for. Pure mathematics is commonly defined as such fields as number theory, algebra, functional analysis, topology/geometry and logic/foundations, as described in the AMS Subject Classification (AMS 2020). On the other hand, applied mathematics is taken to consist of numerical analysis, optimisation, mathematical physics, probability and statistics. That definition cuts across the naïve meanings of “pure” and “applied”, since, for example, number theory can be used in applications like cryptography while methods of calculation in numerical analysis can be studied for their own sake without regard to any uses they may have. An expert on a given topic may investigate it as a pure topic at one point in their career and apply it at another.

From the ethical perspective, it is the naïve meaning of “pure” and “applied” that is relevant. It is, in essence, that we use the labels of \textit{pure} and \textit{applied} mathematics as verbs. On any mathematical topic, one can do pure mathematics (with the aim of developing abstract mathematical truth), and one can do applied mathematics (with the aim of solving an extra-mathematical problem). It is not the mathematical content that makes the work pure or applied, but what it is being done for. 

Pure mathematical research is not free of ethical issues, but they are certainly fewer and less urgent, because its consequences, though occasionally large, are generally more remote and unpredictable. Hersh argued that “in pure mathematics, when restricted just to research and not considering the rest of our professional life, the ethical component is very small” (Hersh 1990, 22). When considering codes of ethics later, we will comment briefly on some issues special to pure mathematics. For now we note that the naïve meaning of “pure” and “applied” does not exempt those trained in the traditionally pure fields from ethical codes, as they must constantly reflect on the focus of their mathematical questions and what their employer and end-users are doing with the mathematical results.

Most calls for a Hippocratic Oath come from mathematicians working on applied mathematics (in the naïve sense), where the impacts are near-term and to some degree foreseeable. Application areas include mathematical modelling of problems in fields like climate, logistics and combat, and a range of mathematical disciplines such as statistics applied to real data, financial and actuarial mathematics, epidemiology, cryptography, operations research, quantitative risk, image compression, mathematical ecology and so on. Since the consequences of mathematical conclusions in these fields are sometimes very large, important ethical issues arise.

\section*{Part I: The evolving landscape of ethics in mathematics}

Early discussions about codes of conduct and oaths for statisticians can already be found during the 1950s. By that time statisticians had been widely employed in a range of applied contexts. From government agencies such as the national census to industrial settings doing quality control, statisticians were constantly challenged in their role. It was only natural that discussions about qualitative statistical work and its purposes would be fairly entangled with ethical questions in the history of statistics. As a consequence, however, the progression of ideas about quality and the nature of the statistical profession have been “nomadic”, jumping between the various sectors where statisticians have been employed (Prévost 2018). The most prominent example of such a discussion occurred in a 1952 special issue of the \textit{American Statistician}, aptly entitled “Standards of Statistical Conduct in Business and Government” (Statistical Conduct 1952). In this volume, William W.K. Freeman suggested a first statistician’s oath, or a set of “principles”:
\begin{quote}
"Like Euclid and all the other great thinkers who have used symbols to reveal the truths of nature, I will be a seeker of the truth. Realising that numbers are only a shorthand convention for describing past events and forecasting trends, I will search for those facts expressed in numbers which show relationships and events most truly. Though surrounded by the clamour of the marketplace or of the political arena, I will not be a fraud, who selects figures to prove by chicanery a misnamed conclusion." (Freeman 1952, 19)
\end{quote}
Freeman saw himself as a seeker of truth in a work environment shaped by economic and political forces. His principles highlight a general theme encountered in discussions about ethics in mathematics. For many mathematicians, the search for mathematical truth ought to withstand and be independent of external pressures because they consider pure mathematics as a neutral and value-free truth. Pure mathematicians can then deduce that stricter codes of conduct matter for the applied fields only. In recent writings on levels of ethical engagement, this stance has been called “level 0: believing there is no ethics in mathematics” (Chiodo and Bursill-Hall 2018, 5). The distinction between the pure and applied becomes even clearer in James Morton’s contribution to the same special issue (Morton 1952, 6-7). Here, he proposed to classify statisticians into three groups:
\begin{enumerate}
    \item The “theorist” who the author did not consider further as “we need not concern ourselves with his conduct”,
    \item the “subject matter specialist” for whom “problems of statistical conduct [...] will seldom become acute”,
    \item and the “fact finding statistician” for whom “the question of standards of conduct becomes important both as a morale builder and as a protective device.”
\end{enumerate}
Only two years later, in 1954, Darrell Huff published his now famous book \textit{How to Lie with Statistics} (Huff 1954), which has sold better than any other statistics book in history (Steele 2005). It put the problems surrounding statistical practice into the public sphere and statisticians continued to grapple with them. The discussions in the American Statistical Association went on until the 1980s (Ellenberg 1983), when it published its first code of conduct (Hurwitz and Gardenier 2012). While constituting a minority by then, some few statisticians still argued against explicitly adopting an ethical commitment. Anderson neglected the societal boundedness of statistical work, reducing moral dilemmas or misconduct to problems of individual ability: “It is often not a matter of statisticians choosing to violate accepted principles of right or wrong. They do so out of ignorance or incompetence” (Anderson 1981, as quoted in Ellenberg 1983). 

However, in the 1990s many other statistical societies had followed and in 1993 the Royal Statistical Society introduced its first code of conduct. It requires its members to “act in the public interest [...] fulfil their obligations to employers and clients [...] to the profession and society [and] at all time show professional competence and integrity” (Royal Statistical Society 1993/2014, 6). Statistics had become the next professionalised area of mathematics, including the ethical commitment characteristic for professions.   

The closely connected field of operations research (OR) traces its origins back to the 1940s. Initially a very applied field, originally in military operations, it developed pure aspects as problems became studied for their own sake. Discussions about ethics in OR go back to at least the 1960s, however, even today most of the discussions surrounding ethics are restricted to the OR process and not the mathematical models (Wenstøp 2010). In practise, ethics in OR is often still only considered in an ad-hoc manner (Ormerod and Ulrich 2013). For an overview of recent developments in OR ethics, we also refer to the work of Le Menestrel and van Wassenhove (2009), but we note here that Brans (2002) suggested a Hippocratic Oath. In his “Oath of Prometheus'' he considered the different roles of OR researcher, teacher, decision-maker and consultant/analyst separately. The oath focuses on a “commitment in favour of fairness, devotion and honesty” and he underscores that “all those taking the oath [...] may be proud of it” (Ibid, 196). Out of it grew the \textit{EURO Working Group on Ethics and OR} (Brans and Gallo 2004), yet by 2009 only 75 researchers had signed the oath (Gass 2009). This lends weight to an observation we discuss later: an oath alone might not be sufficient.

But what were the other mathematical sciences doing during this period? At the end of the Cold War, a small group of mathematicians and scientists started calling for a code of conduct and ethics courses for mathematicians. At a conference in 1988 the American mathematician Chandler Davis asked whether the discipline needed a Hippocratic Oath (Davis 1989). Shortly after, Reuben Hersh wrote \textit{Mathematics and Ethics}, a four-page paper outlining some thoughts about the ethical responsibility of mathematicians (Hersh 1990). He is representative for his colleagues, who often worked alone on such issues and published a few single digit page papers. Their ideas did not find much acceptance in the wider mathematical community of the time. Most striking about his paper is his need to defend himself for talking about ethics in the first place. He began by assuring his readership that he was neither going to tell them what to do, nor that he “solved any big problem regarding mathematics and ethics”, but that he had merely “thought about the question, and in the process of thinking about it has had some ideas which he’d like to offer” (Ibid, 20). That it was still rare even in the 1990s to think about ethics in mathematics becomes abundantly clear from James Franklin’s account about developing such a course for his department at the University of New South Wales (imposed by university policy): 

\begin{quote}
“I was approached by a sheepish Head of School with a message along these lines, ‘We’re not desperate to find someone to create Professional Issues and Ethics in Mathematics; but if you don’t do it, we will be.’ I accepted.” (Franklin 2005, 98)
\end{quote}

Despite the fact that by the 1990s statisticians and the actuarial profession already had well defined codes of conduct and active discussions, and that cryptographers were fighting the US government and its allies in what has long been fittingly called the “Crypto Wars” on restrictions about encryption and cryptographic research, the broader mathematical community remained rather silent on the wider issues of ethics in mathematics. There were some exceptions, where individual mathematicians wrote short pieces on ethical issues arising in mathematical work and how to teach about these, such as Bonnie Shulman's paper \textit{Is There Enough Poison Gas to Kill the City?: The Teaching of Ethics in Mathematics Classes} (Shulman 2002). But most mathematical societies focused their efforts on producing codes of conduct focusing on academic misconduct such as plagiarism and teaching responsibilities.

In the following years, however, individual subfields of mathematics came under pressure and experienced moments of crises. These fields would go on to develop their own forms of Hippocratic Oaths and codes of conduct. After the global financial crisis of 2007-2008, modern financial mathematics had been put into the spotlight. As a consequence, Emanuel Derman and Paul Wilmott published the Financial Modeler’s Manifesto (Wilmott and Derman 2009). It contains the following Hippocratic Oath:

\begin{quote}
“I will remember that I didn’t make the world, and it doesn’t satisfy my equations. Though I will use models boldly to estimate value, I will not be overly impressed by mathematics. I will never sacrifice reality for elegance without explaining why I have done so. Nor will I give the people who use my model false comfort about its accuracy. Instead, I will make explicit its assumptions and oversights. I understand that my work may have enormous effects on society and the economy, many of them beyond my comprehension.”
\end{quote}

After the Snowden Leaks in 2013, Phillip Rogaway (2015) published \textit{The Moral Character of Cryptographic Work}. At the same time, numerous newspaper articles and comments by mathematicians emerged about being employed by intelligence agencies and their work (e.g.~Beilinson 2013, Leinster 2014, Wertheimer 2015, Korman and Tong 2016).  By this time, ethical challenges had appeared in at least six mathematical disciplines: actuarial science, mathematical physics, statistics, operations research, financial mathematics, and cryptography. A seventh was about to follow with the crisis of algorithms and big data. As the impact of social media and the power of big data became evident, the research into its ethics grew rapidly. In her 2016 book \textit{Weapons of Math Destruction}, Cathy O’Neil (2016) laid out many of the dangers surrounding the use of mathematical technologies. She also argued that data scientists, many of whom have come straight out of a mathematics degree with very little additional training, ought to have a Hippocratic Oath, even though “a Hippocratic oath alone is insufficient for the task that lies ahead, because at the end of the day data scientists are not corporations” (Upchurch 2018). Since then many codes of conduct and guidelines have been developed in the fields of data science and artificial intelligence (Jobin, Ienca, and Vayena 2019). 

In 2016, the Cambridge University Ethics in Mathematics Society and, shortly thereafter, the Cambridge University Ethics in Mathematics Project (2018) were founded. Focused on “teaching responsible behaviour and ethical awareness to mathematicians”, they hosted the first two international conferences on ethics in mathematics (“EiM 1” 2018; “EiM 2” 2019). Other research efforts by mathematicians have followed, putting the focus on the issue of mathematics and social justice (Buell and Shulman 2019), or on the values associated with mathematics (Ernest 2016a, 2016b). However, despite a growing awareness, the idea of ethics in all of mathematics has not yet been widely accepted (Chiodo and Clifton 2019). In many places ethical considerations come second and at times are even considered a burden. Hence, the calls for a Hippocratic Oath for mathematicians persist. The latest goes back to Hannah Fry in 2019: 
\begin{quote}
"We need a Hippocratic oath in the same way it exists for medicine. In medicine, you learn about ethics from day one. In mathematics, it’s a bolt-on at best. It has to be there from day one and at the forefront of your mind in every step you take." (Sample 2019)
\end{quote}
It was argued by Chiodo and Bursill-Hall that such teaching of ethics to mathematicians must be woven in throughout a standard mathematics curriculum, needing to be present not only at all times, but in all individual courses in some way, to normalise ethical considerations and make them an everyday part of mathematical work (Chiodo and Bursill-Hall 2019). In the fields of education and philosophy much has been written about the ethics of teaching mathematics and its role at the high school level (e.g.~Dubbs 2020), but it remains a desideratum for the mathematical curricula at universities. 

As we have seen, ethics in mathematics has been interpreted as “ethics in X” for the longest time where X could be an applied subfield encountering ethical challenges. Few publications have thought about the general issue of ethics in higher mathematics. This relates to the fact that mathematics as a whole, in contrast to some of its subfields, does not count as a profession. Mathematics is widely applied in industrial settings (Gr\"otschel, Lucas, and Mehrmann 2010) and mathematicians take on various job titles throughout their career, including software engineer, statistician, or analyst (SIAM 2017).  From the earliest surviving copies of the ancient Hippocratic Oath to today’s widely used version of the World’s Medical Association’s \textit{Declaration of Geneva}, however, all oaths had in common that they understood physicians as a profession. What they did and are doing has always been clearly delineated. What can really be said about mathematicians as a group? Would it make sense to have a Hippocratic Oath for Mathematicians?

\section*{Part II:  A Hippocratic Oath for mathematicians?}
As the starting point for our analysis, we will use the Hippocratic Oath as re-written by Lasagna in 1964 (Various Physicians Oaths 2016), which is currently used by a third of US medical schools (Crawshaw et al. 2016). Lasagna’s oath is more detailed than the more widely applied version of the World Medical Association’s \textit{Declaration of Geneva} (WMA 2021), and thus provides a suitable case study.  In response to Brad Smith and Harry Shum (Microsoft 2018), Etzioni (2018) already edited it for the case of artificial intelligence. However, since mathematics is a broader discipline, encompassing much of artificial intelligence, we will not use Etzioni’s rewritten version as the starting point of our discussion. Otherwise, we are prone to overlook the subtleties special to mathematics and its wider community. Instead, we will go through each of the ten sections, and analyse if and how Lasagna’s oath could be adapted to the specific needs of mathematicians. 
\newline\newline
\textbf{1. I swear to fulfil, to the best of my ability and judgement, this covenant.}
\newline\newline
At the beginning of their career, and by taking this oath, physicians enter a contractual agreement (“covenant”) with their peers and their national accreditation body. Neither party to this contract exists in the same form for mathematics: most mathematical societies do not (yet) give out the accreditation of being a “chartered mathematician” and the mathematical community is comprised of many sub-communities, some of which do not even share the same degree title.     
\newline\newline
\textbf{2. I will respect the hard-won scientific gains of those physicians in whose steps I walk, and gladly share such knowledge as is mine with those who are to follow.}
\newline\newline
Almost all mathematical work builds on the work of others. Thus, this part of the oath could appeal to mathematicians not to exploit the mathematical works and developments of others and misuse them in ways that cause harm. Such concerns have already been raised by the mathematician Andrew Wiles, regarding the misuse of mathematical tools leading up to the 2008 financial crisis (Devlin 2013).

Mathematicians ought to share, publish and disseminate their mathematical knowledge widely. However, given how much of our modern society is dependent on mathematical technologies, this instruction might include what computer scientists know as “responsible disclosure”. Not all cryptographic results can be immediately published without endangering the existing critical infrastructure. Responsible disclosure has also been suggested for artificial intelligence (Hagendorff 2020). The issue of what constitutes the mathematical community still stands.
\newline\newline
\textbf{3. I will apply, for the benefit of the sick, all measures [that] are required, avoiding those twin traps of overtreatment and therapeutic nihilism.}
\newline\newline
Despite an increasingly holistic and proactive look at human health, most medical work is still concerned with trying to heal the sick. The tools of mathematics, however, are often used to improve already working solutions and are not necessarily about restoration. In addition, mathematicians might take inspiration from the Financial Modeler’s manifesto and interpret (3) as an instruction to use the “right” amount, and type of, mathematics: mathematicians should neither over-model or over-solve a problem, nor give up too easily.
\newline\newline
\textbf{4. I will remember that there is art to medicine as well as science, and that warmth, sympathy, and understanding may outweigh the surgeon's knife or the chemist's drug.}
\newline\newline
This could be a reminder to mathematicians to remember their humanity, and that their work should not be applied strictly and rigidly at all times. Applied mathematical work entails non-mathematical decisions, asking questions and understanding when mathematics is not the right tool to responsibly solve the problem. In an age of quantified decision making, mathematicians are still faced with imprecise assumptions to make and steps to take, which require perspective, experience, and humanity. Those affected by mathematical technologies should not be “reduced to a number.”
\newline\newline
\textbf{5. I will not be ashamed to say "I know not," nor will I fail to call in my colleagues when the skills of another are needed for a patient's recovery.}
\newline\newline
Doctors are trained to seek the opinion of other doctors when they feel their medical knowledge is lacking. This is similarly true for mathematicians. However, most extra-mathematical problems, situations and scenarios require extensive domain knowledge to be resolved. There may be times when it is necessary to seek help, guidance, and insight not just from other mathematicians, but also from experts in other fields (Bennett Moses 2018). At the beginning of the Covid-19 pandemic, mathematicians tried to model the transmissions and spread of Covid in elderly care homes. However, it has been shown that the modellers underestimated how poor the shielding was between care homes and the wider society, and that they did not ask about it during the modelling phase. No representatives from the care homes were asked for specific advice, and consequently the models missed the impending outbreaks throughout the care facilities (BBC Two 2020).  The interaction with domain experts is of vital importance to succeed and guidelines for responsible modelling should be followed (Chiodo and Müller 2020). 
\newline\newline
\textbf{6.a. I will respect the privacy of my patients, for their problems are not disclosed to me that the world may know. Most especially must I tread with care in matters of life and death.}
\newline\newline
Mathematicians sometimes deal with sensitive issues and confidential data. Mathematical tools such as alternative credit scoring, for instance, use vast amounts of data and in return can severely impact people’s lives, their decision making, and the choices available to them. To fully respect their autonomy demands that developers respect data privacy and understand that the effects of such data still linger in the trained algorithm long after the initial data has been deleted. It also requires developers to understand more complex issues in data justice (for an overview, see Taylor 2017 and Dencik et al. 2019).  Mathematicians need to to be aware of how the decision rules, patterns and associations inferred from a data set pertaining to a subset of individuals can still be used to influence, and possibly exploit, a larger group, and thus could be a violation of both individual and group privacy (Floridi 2014, Suh et al. 2018).

It is not uncommon to have initially studied mathematics out of enjoyment, and to see it as a game or even as a competitive sport. Consider, for example, the International Mathematical Olympiad which has been running annually since 1959, and its many national counterparts. Examples such as these, however, show that dealing with sensitive data can have real-life consequences and substantial impact which the mathematician should recognise.
\newline\newline
\textbf{
6.b. If it is given to me to save a life, all thanks. But it may also be within my power to take a life; this awesome responsibility must be faced with great humbleness and awareness of my own frailty. }
\newline\newline
Mathematics has the potential to do vast good, but also to do vast harm. However, as is often the case, for mathematics it is less clear than in medicine what constitutes good and bad. Mathematicians often work in complex environments and are faced with competing interests, market and political forces. 
\newline\newline
\textbf{
6.c. Above all, I must not play at God.}
\newline\newline
Predictive algorithms, pricing systems, insurance and many other mathematical systems can impact millions of people. The potential to “play at God” is ever present. It should be understood and avoided.
\newline\newline
\textbf{7. I will remember that I do not treat a fever chart, a cancerous growth, but a sick human being, whose illness may affect the person's family and economic stability. […]}
\newline\newline
When mathematicians are asked to solve a societal or industrial problem, they should acknowledge that their mathematics is a means to an end, and not an end in itself. Their work can have an impact on people’s private lives, their economic status and society more generally. Broader considerations such as bias, fairness, rights and values must be in the mathematician’s mind, for their work not to become “algorithms of oppression” (Noble 2018). However, this point also raises a crucial subtlety which we will discuss later: who is a mathematician’s “patient”? How does a fiduciary duty, and a duty of care, look for mathematicians, and to whom do they owe it?
\newline\newline
\textbf{
8. I will prevent disease whenever I can, for prevention is preferable to cure.}
\newline\newline
This could be viewed as mathematicians acknowledging that preventing problematic outcomes from their work, in advance, is always preferable to Silicon Valley’s mantra of “moving fast and breaking things” (Levy 2014).
\newline\newline
\textbf{
9. I will remember that I remain a member of society, with special obligations to all my fellow human beings, those sound of mind and body as well as the infirm.}
\newline\newline
Mathematicians have all the moral obligations to society that everyone else has, and this is not limited to the entity they work for. This has already been included in codes of ethics for many statistical societies.
\newline\newline
\textbf{
10. If I do not violate this oath, may I enjoy life and art, respected while I live and remembered with affection thereafter. […]}
\newline\newline
The oath could conclude by assuring mathematicians that, by following it, they will live a happy and virtuous life, be remembered well, and be free to enjoy doing mathematics and seeing its positive impact. Mathematicians are very proud of their discipline, demonstrated in the way they name theorems and results. Any mathematician acting in line with the oath can know that they are leaving a positive legacy. Mathematicians are able to trace their educational roots back over 1000 years (Mathematics Genealogy Project 2021), and unlike in other disciplines, results or ideas are almost never deleted or overwritten. Their output might not receive much attention, but it does get “added to the shelf” forever.

\mbox{}

Versions of the Hippocratic Oath for medical practitioners have been in existence for over 2000 years. That has had a deep and lasting impact on the culture of medical practitioners. They know that ethics is important in their field, and they swear an oath which was sworn by their teachers, and previous generations. In mathematics, there is no such oath. But more than that, the idea of ethical responsibility has not permeated into the disciplinary culture, into the training, or into the mentorship, and the phrase “back in my day we had no such oath” may survive for years to come.

Present in the original Oath of Hippocrates is the notion to “do no harm” (Various Physicians Oaths 2016). This is often interpreted as to only prescribe treatments when the benefits outweigh the risks; harming the body through surgery is allowed in this context. How might mathematicians interpret this when it is not entirely clear what harm actually means? It might be easy for mathematicians to restrict this notion to their employer’s interests, but how should they consider impacted parties further afield? Human Rights issues are already being explored in the context of artificial intelligence (Risse 2019), and it must also be put on the agenda for mathematics: How can mathematical work violate human rights, and how can mathematicians take countermeasures?

The Hippocratic Oath lacks detail of how a medical practitioner should act. But this is supplemented by various medical ethics standards, and medical ethics training that is received by almost all medical students. Furthermore, the medical profession is accredited and governed by extensive rules and regulations on what warrants improper conduct. So, in essence, the oath is merely an easy-to-remember reminder of an otherwise deep, complicated and well-supported system. The oath alone, though useful, would probably not achieve all of it. Mathematicians have no such systems. There are remarkably few rules or regulations governing general mathematical work, there is no widely-recognised accreditation process that can be revoked, and there is almost no ethical training offered to mathematicians as part of their degree courses. As even the existence of ethics in mathematics is not yet widely accepted, its content would likely be subjected to extensive challenges.

\section*{Part III: What is the alternative?}
\subsection*{Ethical principles from law and codes of ethics}
In this section we extract from existing legal/moral principles and from codes of ethics in quantitative fields some general principles applicable to the ethics of any mathematical work that impacts the public. The long experience of the legal world and of practical disciplines such as engineering and statistics has resulted in some well-grounded thought on ethical principles. The mathematical community can stand on their shoulders. Irrespective of whether or not a code of ethics for mathematics is desirable, the long-term reflection that is distilled in the law and the codes of cognate disciplines has valuable lessons for mathematicians.  We begin with the most general principles.

The law recognises one very general ethical principle applicable across the board to professional dealings that impact the public – duty of care.    The principle of duty of care is a simple one. A manufacturer or provider of services who does something that can be reasonably foreseen to cause harm is liable for that harm. That is so even if the harm is “downstream”, rather than to the party with whom the manufacturer or provider contracted. In one of the first cases to establish the principle in American law, \textit{MacPherson v.~Buick Motor Co.~}(1916), a driver had been injured by the collapse of the wooden wheel of his Buick, bought from a dealer. In holding the manufacturer liable, Judge Cardozo said “If to the element of danger there is added knowledge that the thing will be used by persons other than the purchaser, and used without new tests, then, irrespective of contract, the manufacturer of this thing of danger is under a duty to make it carefully.” In the classic British case of the “snail in the bottle,” \textit{Donoghue v.~Stevenson} (1932), where Mrs Donoghue became ill through ingesting a decomposed snail in the bottle of ginger beer manufactured by Stevenson but sold in a café, Lord Atkin held:

\begin{quote}
"You must take reasonable care to avoid acts or omissions which you can reasonably foresee would be likely to injure your neighbour. Who, then, in law, is my neighbour? The answer seems to be – persons who are so closely and directly affected by my act that I ought reasonably to have them in contemplation as being so affected when I am directing my mind to the acts or omissions which are called in question." (Donoghue v Stevenson (1932): 44)
\end{quote}
While the law of negligence has since become complex on questions such as the degree of remoteness of causality and of foreseeability, and also differs somewhat between jurisdictions, the basic principle remains clear. The ethical requirement to avoid reasonably foreseeable harm caused by the products or services one provides translates into legal liability.

That applies directly to the products supplied by mathematicians that could have detrimental effects on the public. In particular, a failure to make an effort to foresee effects is not a legally sufficient excuse for not foreseeing effects that are, as a matter of objective fact, reasonably foreseeable.

While the language of “duty of care” is not found in codes of ethics of the quantitative professions, a near-equivalent idea appears in some codes under the name “the public interest”.  That concept refers to both harmful and beneficial effects of work on the wider public.  For example, the Royal Statistical Society puts first among its principles:

\begin{quote}
\begin{itemize}
\item "Fellows should always be aware of their overriding responsibility to the public good; including public health, safety and environment."
\item "A Fellow’s obligations to employers, clients and the profession can never override this [...]" 
\item "Fellows shall in their professional practice have regard to basic human rights and shall avoid any actions that adversely affect such rights." (Royal Statistical Society 1993/2014)
\end{itemize}
\end{quote}

The Australian Computer Society Code (2014) also starts with the primacy of the public interest, which it fills out with: 
\begin{quote}
\begin{itemize}
\item "identify those potentially impacted by your work and explicitly consider their interests" (similar in Institute of Mathematics and its Applications 2018).
\end{itemize}
\end{quote} 
Engineers Australia (2019) advises:
\begin{quote}
\begin{itemize}
\item "Practise engineering to foster the health, safety and wellbeing of the community and the environment [...] incorporate social, cultural, health, safety, environmental and economic considerations into the engineering task. [...] Balance the needs of the present with the needs of future generations"
\end{itemize}
\end{quote}

Those comments on the public interest are more positive in tone than the legal concept of duty of care, which is more narrowly concerned with avoidance of harm. Nevertheless, both concepts relate to the ethical requirement to take into account significantly and primarily the interests of all “stakeholders” who are reasonably foreseen to be affected by work, not just the interests of the contracting parties who commission and perform the work.

Before considering some other principles found in codes, it is useful to distinguish between \textit{codes of ethics} and \textit{codes of conduct}, even though that distinction is not clearly made in the names of existing professional codes. A code of ethics is the broader concept. It includes principles applicable to all the ethical aspects of a professional’s work, including its relation to the public interest. A code of conduct, on the other hand, restricts itself to the personal behaviour of a professional in the field, enjoining such virtues as honesty, competence, respect between professionals and ensuring work does not bring the profession into disrepute. In this terminology, the Royal Statistical Society’s code, though called a code of conduct, is broad enough to be a code of ethics. But the codes of the American Mathematical Society, the Institute of Mathematics and its Applications, and the Society of Industrial and Applied Mathematics, are codes of conduct.

These issues of professional conduct are secondary to those involving the public interest, since there is no ethical point in behaving with honesty, integrity and competence on a project that is itself harmful. Nevertheless, they are crucial for working ethically in a profession, and several codes have expressed requirements such as honesty/integrity, competence, avoiding conflict of interest, engaging/communicating, and attending to data issues:
\\\mbox{}
\\\textbf{Honesty/integrity:} Professionals should state truly the facts, about such matters as what is revealed by their work, their own competence and qualifications, foreseen impacts of the work, correct attribution of work, and the likelihood of contracted work being able to reach credible conclusions. As the SIAM code puts it, “act with integrity, and strive to be objective, unbiased, and truthful in all aspects of our work.” (SIAM n.d.)
\\\mbox{}
\\\textbf{Competence/diligence:} The work of any professional should be performed with competence and diligence. “Fellows shall carry out work with due care and diligence in accordance with the requirements of the employer or client,” says the Royal Statistical Society (1993/2014, 2). But high standards are expected especially of mathematicians. Mathematics holds out the prospect of absolutely certain, proven results. The public trusts that. Where results are provable, they should be proved, and where that is impossible, they should be to the highest standards. Competence includes updating one’s professional knowledge with new developments, consulting domain experts and knowing the applicable law.
\\\mbox{}
\\\textbf{Conflict of interest:} Where a professional has or is reasonably perceived to have a conflict between the interests of the client or public and those of him/herself or his “kin” (such as family, students or colleagues), the conflict must be appropriately managed. That could mean appropriate disclosure or recusing oneself from decision-making. If the conflict is severe, it could mean refusing work. The Australian Computer Society advises:
\begin{quote}
\begin{itemize}
\item"[to] raise with stakeholders any potential conflicts between your professional activity and legal or other accepted public requirements"
\item "[to] advise your stakeholders as soon as possible of any conflicts of interest or conscientious objections that you have" (ACS 2014)
\end{itemize}
\end{quote}
The Actuaries Code states:
\begin{quote}
\begin{itemize}
\item "Members must ensure that their professional judgement is not compromised, and cannot reasonably be seen to be compromised, by bias, conflict of interest, or the undue influence of others."
\item "Members must take reasonable steps to ensure that they are aware of any relevant interests that might create a conflict."
\item "Members must not act where there is an unreconciled conflict of interest. (Institute and Faculty of Actuaries 2019)"
\end{itemize}
\end{quote}
\mbox{}
\\\textbf{Engaging/communicating:} Professional work is a shared enterprise, so it requires engaging and communicating with all affected stakeholders. The client has a right and duty to communicate fully what it wants done and the relevant information for doing it. Other stakeholders such as those affected by the implementation of an algorithm should be surveyed as well as having their interests guessed at. The nature of multidisciplinary work means that mathematical skills by themselves are rarely sufficient for a task, so mathematicians must recognise where their expertise ends and seek out those with the necessary complementary skills. For example, an algorithm to be implemented in software for public use will need expertise in human-computer interaction to ensure its usability and safety. 

Information on results should be disseminated as freely as possible, subject to any necessary requirements of confidentiality. The American Mathematical Society states

\begin{quote}
\begin{itemize}
    \item “Freedom to publish must sometimes yield to security concerns, but mathematicians should resist excessive secrecy demands whether by government or private institutions.” (AMS 2005/2019) 
    \end{itemize}
\end{quote}

For example, the results of an aborted or suppressed statistical trial which proves inconvenient to a client’s interests might need to be disclosed if they indicate a danger to the public. (McGoey and Jackson 2009)
\\\mbox{}
\\\textbf{Data issues:} Much of applied mathematical work is based on or at least takes into account quantitative data, for which special issues arise. The data itself must be ethically collected, with consent where appropriate. Mathematicians must be aware of, and advise clients and the public of, the sufficiency of the data for the claims based on it and the need to guard against biases in it. The Oxford-Munich Code for Professional Data Scientists (Grindrod and Moreno 2018) lists a number of special issues arising in making inferences from data reliably and ethically, including:
\begin{quote}
\begin{itemize}
\item "The Data Scientist is responsible for assessing the adequacy of data to solve the particular problem and to share the results of the analysis, indicating any risks or potential implications due to lack of data quality or availability."
\item "The Data Scientist is responsible for separating correlations that are the results of chance or deliberate data-mining driven searches vs. well established hypothesis-driven correlated information."
\item "The Data Scientist shall retain copies of the original data unaltered."
\item "The Data Scientist shall not apply any technique (combination, enriching, etc) to turn information that has been designed to be “de-identifiable” into “identifiable” again"
\end{itemize}
\end{quote}
While these are in the first instances technical issues about inference from data, it is the expert knowledge of them held by mathematicians, statisticians and data scientists that gives them a cognitive advantage over clients and users, and hence a responsibility to act competently and honestly.
\subsection*{A Code of Ethics for Mathematicians?}
It would certainly be possible to construct a code of ethics for mathematicians, adapting recognised principles from the existing codes of ethics, Hippocratic oaths and adding some advice on detailed issues that are special to mathematics. Is that a useful exercise?

It is true that such a code would not (yet) be enforceable. Mathematics is not, and is unlikely to soon become, fully professionalised like medicine, law and the actuarial profession, with accreditation and a powerful professional association capable of disciplining members. 

The institutionalised actuarial community traces its origins back to the nineteenth century. In the UK the origins of the Institute and Faculty of Actuaries (IFoA) goes back to a Royal Charter in 1848 (Privy Council 2021), and it still works as the only institution to both represent and regulate actuaries in the UK (IFoA n.d.). Their work surrounding accreditation and regulation is highly effective. Similarly, statisticians have undergone a professionalisation with accreditation mechanisms. The Royal Statistical Society has existed with Royal Charter in the UK since 1887 (Privy Council 2021). It accredits and governs the actions of professional statisticians in the UK following an extensive code of conduct, and their work further combats various misuses of statistics in society, such as in the case of Sally Clark (Nobles and Schiff 2005). But as we have seen, for mathematicians the situation is much more complicated. Even the labelling of mathematicians as mathematicians is problematic, with many graduates in mathematics later identifying as software engineers or computer scientists, sometimes without further professional training. It may therefore be unclear who such a code would apply to.

Furthermore, our analysis of Lasagna’s Hippocratic Oath has introduced two additional problems: existing oaths for physicians rely on there being a well-defined patient, and a common understanding of “good” and “harm”. In medicine, as in the other classical professions like law and the clergy, a well-defined client has a relationship of personal encounter with the professional, who has explicitly agreed to act in the client’s interest. (Koehn, 1994, ch. 1, Franklin, 2011); that personal relationship structures the professional ethics of the relationship. For mathematicians (as for engineers and computer scientists), we cannot simply define a patient, but have to consider a variety of stakeholders in addition to the party employing the mathematician: those using the technology (e.g.~the users of a search engine), the customers of the company (e.g.~companies buying ad space on social media), their employer and the various funding sources, and other impacted parties (e.g.~those being advertised/sold to). 

Mathematicians are constantly confronted with trade-offs between the stakeholders and within the individual groups. Recent research has shown that even considerations such as fairness and bias are often interpreted with different and competing notions and results: individual fairness does not necessarily imply or contradict group fairness, nor vice versa (Binns 2020). Mathematics amplifies the reach of technologies and the situation complexifies quickly, and those impacted by it can lie in the thousands if not millions. Here lessons could be taken from existing codes for statisticians and the Oxford-Munich Code for Data Science as both communities deal with similar issues.

A further difficulty arises in the proximity and abstraction of impacted parties. Unlike mathematicians, a physician can usually see the patient and inquire with them to gain a better understanding. Mathematical work usually affects many individuals and groups far away from and unseen by the mathematician, and often in ways that the mathematician might never see nor experience first hand. Complex mathematical systems lack feedback mechanisms, of the kind that exist in a doctor-patient relationship, and are difficult to change if mistakes were made. Questions of consent and accountability are already highly debated issues in artificial intelligence and the wider field of ethics of algorithms. All their challenges remain: Can users object and withdraw consent in the same way a physician’s patient could? Do alternative technologies exist and what are the opportunity costs associated with them? Who should be held accountable if things go wrong, and what are the best ways to establish accountability? 

Codes and oaths partially take their value from being a source of reference, for example in litigation. Computer scientists, who face the same problem of a non-professionalised and diversified discipline, have found some value in a code. The Australian Computer Society comments, under “relevance to law”: “This Code of Professional Conduct has relevance to professional standards legislation. Failure to abide by the Code could be used as grounds for a claim of professional negligence. The Code may be quoted by an expert witness giving an assessment of professional conduct.” (ACS 2014, 4) So one point of a code is that no-one can then say in court “Well, no-one \textit{told} me I was supposed to be ethical, or what that would mean.” If the law inquires whether a mathematician has breached a duty of care by causing foreseeable harm, a code could add flesh to that general principle. For it to be effective, however, requires its authors to engage with the previously mentioned challenges.

Furthermore, the notions of “good” and “harm” are muddied for many mathematical applications. While a clinical patient is usually fighting a disease or ailment and the physician assists in that fight, mathematicians are not necessarily fighting something generally agreed upon as “bad”. Often, problems and moral questions in mathematical work only become obvious after implementation. Few had thought about the role of mathematics in targeted political advertising before the case of Cambridge Analytica, in which the personal data of millions of Facebook users were harvested without their consent (Chan 2019). Arguably, much of modern mathematical research has worked for the benefit of society, but the problems practitioners deal with are so complex that when challenges do occur, the stakes are high. Then mathematicians have to weigh up different notions of “good” and “harm”: public, individual, environmental, or corporate – as do the already existing codes studied above.

Finally, any successful code of ethics must also address the issues of pure mathematical research, which is more comparable to philosophy and art than to engineering. Both philosophy and art bring their own set of ethical challenges, which we have not further considered in this paper. We will, however, briefly need to mention the opportunity cost of pure mathematics and talk about mathematical education before moving on. Is it justified to have such a proportion of very talented and highly trained people devoting their lives to proving theorems when they could be doing something more ethically positive, like medical research? (Franklin 1991) Others argue that pure mathematics must be able to justify itself as a prime expression of the human spirit or essential to human flourishing. (Su 2020). Beyond that, ethical issues arise from the domination of mathematical pedagogy by pure mathematics and the usual internal issues for any profession, such as fairness in hiring and correct attribution of work – issues laid out in the American Mathematical Society’s Code of Ethics. (AMS 2005/2019) It remains an open question if and how a code of ethics for mathematicians should address these issues.

Reaching consensus on a code of ethics for all mathematicians could take a long time, and an effective middle-ground could be institutional mission statements. It has been proposed that individual mathematicians come up with personal mission statements (Harris and Vega 2020), but this approach might be too limited in reach and should be expanded. Individual universities, departments and other employers might introduce or expand their mission statements to go beyond teaching, inclusivity and diversity. These can be tailored to the particular circumstances of the institution and its surroundings, taking into account the relevant cultural and social needs. Such statements might not be universally applicable, but they are quicker to develop, implement, and evolve. Furthermore, they could be tied to existing institutional structures, making it easier to implement tools of accountability and enforcement mechanisms. But most importantly, they would address a concern commonly heard in the context of ethics in mathematics: “Whose ethics?” Those of the people within the institution. With time, mission statements could feed into and provide inspiration for codes of ethics for national and international mathematical societies.

Fleischmann, Hui and Wallace (2017) conducted interviews with computational modellers about human values and codes of ethics. The interviews suggest that a bottom-up approach for codes of ethics could be most promising, as it seems to be supported by modellers and has the potential to include multiple stakeholders and regional values. To our understanding, no similar study has been undertaken for mathematicians, but we suggest that a bottom-up approach using mission statements would be effective, and could provide a suitable space for minorities to be heard. For mission statements to be successful, they need to build on a set of informed, diverse and representative contributors, including those currently institutionally underrepresented.

\section*{Conclusion and recommendations}

We have seen that the development of ethics in the mathematics community remains in a primitive state, at the same time as the impact of mathematical technologies on society has reached an all-time high. We have also seen that there is much to learn from how cognate quantitative disciplines have approached their ethical responsibilities, as well as from general legal and moral principles. In the light of that, we make several recommendations, ordered not by their perceived importance, but instead by how we believe they might be most effectively and rapidly implemented.

Firstly, individual institutions, such as universities and companies employing mathematicians, should draw up or modify their own small-scale mission statements to specifically include references to ethical considerations of the mathematical work relating to their activities. These can be done quickly, and can evolve in time.

Secondly, the professional societies should encourage both the inclusion of standalone courses on professional issues and ethics in all accredited mathematics degrees, and the incorporation of ethical awareness and training into all mathematics teaching, courses and modules. They should sponsor the development of syllabuses for such courses, with internationally common core content but with variations for local conditions and allowing for advice from industry practitioners.

Thirdly, the mathematical societies, in both pure and applied mathematics, should collaborate to sponsor work towards an international code of ethics for working mathematicians. It should go beyond a mere code of conduct for the internal workings of the academic mathematics profession, but instead foreground the primacy of the public interest. It should distinguish between the ethics applicable to pure and applied work. The code when finalised should become an official policy document of each society.

And finally, the professional societies should institute complaints procedures to which allegations of unethical professional behaviour by mathematicians can be referred for investigation. If allegations are found to be proved, an official statement by the society should make public the results of the investigation. Despite not carrying any legal standing, they may still deter unethical practise among mathematicians through general community awareness. While the development of mission statements, courses, and codes of ethics is crucial, it is eventually needed to go beyond passive documents, and to create an environment actively conducive to ethical mathematical practise.

We finish by acknowledging that some in the mathematical community might be tempted to let the perfect be the enemy of the good, while they are waiting for a perfect universal solution to be implemented. However, choosing to do nothing is still a choice, and is probably far from an ideal solution.

\nocite{*}
\bibliographystyle{apacite}
\bibliography{allrefs.bib}

\end{document}